\input amstex
\documentstyle{amsppt}

\magnification=\magstep1

\pageheight{9.0truein}
\pagewidth{6.5truein}
\NoBlackBoxes
\TagsAsMath

\input xy \xyoption{matrix} \xyoption{arrow} \xyoption{curve} 

\def\NN{{\Bbb N}}
\def\ZZ{{\Bbb Z}}

\def\E{{\Cal E}}

\def\calP{{\Cal P}}

\def\Ehat{\widehat{E}}
\def\Fhat{\widehat{F}}
\def\etahat{\widehat{\eta}}
\def\phihat{\widehat{\phi}}
\def\calEhat{\widehat{\E}}

\def\pibar{\overline{p_i}}
\def\qbar{\overline{q}}

\def\vbar{\overline{v}}
\def\wbar{\overline{w}}

\def\pjstarbar{\overline{p^*_j}}
\def\qjstarbar{\overline{q^*_j}}

\def\xbar{\overline{x}}

\def\DiGr{\operatorname{\bold{DiGr}}}
\def\KAlg{K\text{-}{\operatorname{\bold{Alg}}}}
\def\CKGr{\operatorname{\bold{CKGr}}}
\def\Fld{\operatorname{\bold{Fld}}}
\def\Rng{\operatorname{\bold{Rng}}}
\def\length{\operatorname{length}}
\def\card{\operatorname{card}}
\def\End{\operatorname{End}}
\def\bull{$\bullet$}
\def\dirlim{{\underset\longrightarrow\to\lim}\,}
\def\RMod{R{\operatorname{{-}\bold{Mod}}}}
\def\IMod{I{\operatorname{{-}\bold{Mod}}}}
\def\JMod{J{\operatorname{{-}\bold{Mod}}}}
\def\IJMod{(I/J){\operatorname{{-}\bold{Mod}}}}
\def\FP{\operatorname{\bold{FP}}}
\def\id{\operatorname{id}}

\def\AAa{{\bf 1}}
\def\AAb{{\bf 2}}
\def\AAc{{\bf 3}}
\def\AAPS{{\bf 4}}
\def\AAS{{\bf 5}}
\def\Ara{{\bf 6}}
\def\AraB{{\bf 7}}
\def\AGS{{\bf 8}}
\def\AGOP{{\bf 9}}
\def\AMP{{\bf 10}}
\def\AP{{\bf 11}}
\def\AMMS{{\bf 12}}
\def\APS{{\bf 13}}
\def\Malaga{{\bf 14}}
\def\Cun{{\bf 15}}
\def\DT{{\bf 16}}
\def\Fai{{\bf 17}}
\def\GS{{\bf 18}}
\def\GW{{\bf 19}}
\def\Lam{{\bf 20}}
\def\Lea{{\bf 21}}
\def\Nic{{\bf 22}}
\def\Sil{{\bf 23}}
\def\Tom{{\bf 24}}
\def\War{{\bf 25}}

\topmatter

\title Leavitt Path Algebras and Direct Limits \endtitle

\author K. R. Goodearl \endauthor

\address Department of Mathematics, University of California, Santa
Barbara, CA 93106 \endaddress

\email goodearl\@math.ucsb.edu \endemail

\abstract An introduction to Leavitt path algebras of arbitrary directed
graphs is presented, and direct limit techniques are developed, with
which many results that had previously been proved for countable graphs
can be extended to uncountable ones. Such results include
characterizations of simplicity, characterizations of the exchange
property, and cancellation conditions for the K-theoretic monoid of
equivalence classes of idempotent matrices.
\endabstract

\subjclassyear{2000}
\subjclass 16S99, 16G20 \endsubjclass

\endtopmatter

\document

\head Introduction\endhead

The algebras of the title descend from algebras constructed by W. G.
Leavitt \cite{\Lea} to exhibit rings in which free modules of specific
different finite ranks are isomorphic. Later, and independently, J. Cuntz
introduced an analogous class of C*-algebras \cite{\Cun}. Generalizations
of these led to a large class of C*-algebras built from directed graphs,
and the construction was carried to the algebraic category by G. Abrams
and G. Aranda Pino \cite{\AAa}. For both historical and technical
reasons, the graphs used in constructing graph C*-algebras and Leavitt
path algebras have been assumed to be countable, although the
construction does not require this. Our motivation for this paper was to
initiate the study of Leavitt path algebras of uncountable graphs and to
show that many (perhaps all) of the known theorems hold for uncountable as
well as countable graphs. 

The first section of the paper is an expository account of the basic
ideas and constructions involved with Leavitt path algebras $L_K(E)$,
where $K$ is a field and $E$ a (directed) graph. Further expository
material is incorporated into the second section, where appropriate
categories are defined in which direct limits exist. In particular, we
show that $L_K(E)$ is a direct limit of Leavitt path
algebras of certain countable subgraphs of $E$, over a countably directed
set. This result allows ``finitely definable'' properties to be
transferred from the countable to the uncountable case. Section 3
illustrates this procedure with various results concerning ideals. For
instance, we extend a theorem of Tomforde \cite{\Tom} saying that the
ideals of $L_K(E)$ are homogeneous (with respect to the canonical
grading) if and only if $E$ satisfies the graph-theoretical condition
(K), and a theorem of Abrams-Aranda \cite{\AAc} and Tomforde \cite{op\.
cit.} characterizing simplicity of $L_K(E)$ in terms of graph-theoretical
conditions on $E$. In Section 4, we extend another theorem of
Abrams-Aranda \cite{op\. cit.}, showing that $L_K(E)$ is an exchange ring
if and only if $E$ satisfies condition (K). The final section addresses
the abelian monoid
$V(R)$ associated with any ring $R$, which can be built from either
equivalence classes of idempotent matrices or isomorphism classes of
finitely generated projective modules. We give an expository account of
$V(R)$ in the general (nonunital) case, and develop some nonunital module
category machinery in order to show that Morita equivalent idempotent
rings have isomorphic $V$s. This is needed in our final result, 
extending a theorem of Ara-Moreno-Pardo \cite{\AMP}, which says that
$V(L_K(E))$ is always an unperforated, separative refinement monoid.

In order to keep this paper to a reasonable size, we have chosen to
present only a sample of results where a countability hypothesis on the
graph can be removed, and we have not addressed potential C*-algebra
analogs of the methods. We invite readers to explore removing
countability assumptions from other results. In particular, some natural
candidates can be found in \cite{\AAb, \AAc, \AMP, \AP, \APS}.

\head 1. Basics \endhead

In this section, we give some basic notation for graphs, path algebras,
and Leavitt path algebras, and discuss a few basic results about these
objects. For more on the historical background, we recommend Abrams'
article in \cite{\Malaga}. 

\definition{1.1\. Graphs} A {\it directed graph\/} is a 4-tuple $E=
(E^0,E^1,r_E,s_E)$ consisting of two disjoint sets $E^0$, $E^1$ and two
maps $r_E,s_E: E^1\rightarrow E^0$. The adjective ``directed'' is often
omitted, as are the subscripts on the maps $r_E$ and $s_E$, unless
several graphs are under discussion at once. The elements of $E^0$ and
$E^1$ are called the {\it vertices\/} and {\it edges\/} of $E$,
respectively. If $e\in E^1$, then $s(e)$ and $r(e)$ are called the {\it
source\/} and {\it range\/} of $e$, respectively. One also says that {\it
$e$ goes from $s(e)$ to $r(e)$\/}, written $e: s(e)\rightarrow r(e)$ or
drawn $s(e) @>{e}>> r(e)$ in diagrams. 

A {\it directed path\/} in $E$, usually just called a {\it path\/}, is a
sequence of edges with the source of each edge matching the range of
an adjacent edge. There are two conventions for orienting such sequences:
right to left or left to right, and unfortunately both are in common use.
Since the left to right convention is used in a majority of papers on
Leavitt path algebras and graph C*-algebras, we follow that convention
here. Thus, a path in $E$ consists of a sequence $e_1$, $e_2$,
\dots, $e_n$ of edges from $E^1$ such that $r(e_i)=s(e_{i+1})$ for all
$i=1,\dots,n-1$. We shall write such a path as a product $p=e_1e_2\cdots
e_n$, labelling $s(p)=s(e_1)$ and $r(p)= r(e_n)$. The {\it length\/} of a
path is the number of edges it contains. It is important to allow a path
of length zero at each vertex $v$; this is a path with source and range
$v$, including no edges, and the natural notation for this path is just
$v$.
\enddefinition

\definition{1.2\. Duals and Doubles}
The {\it dual\/} of a graph $E$ is a graph $E^*$ consisting of the same
vertices as $E$ but with all edges reversed. Specifically,
\roster
\item"\bull" $(E^*)^0 = E^0$;
\item"\bull" $(E^*)^1$ is a set $(E^1)^*= \{e^*\mid e\in E^1\}$ which is
in bijection with $E^1$ via $e\mapsto e^*$;
\item"\bull" $r(e^*)= s(e)$ and $s(e^*)= r(e)$ for all $e\in E^1$.
\endroster 
The $*$ notation is extended from edges to paths in the obvious manner.
Thus, if $p=e_1e_2\cdots e_n$ is a path in $E$ from a vertex $v$ to a
vertex $w$, then $p^*= e_n^* e_{n-1}^* \cdots e_1^*$ is a path in $E^*$
from $w$ to $v$. (Note the length zero case: $v^*=v$ for $v\in E^0$.)

It is assumed (sometimes only tacitly) that the set $(E^1)^*$ is disjoint
from $E^0\sqcup E^1$. Then the {\it double\/} (or {\it extended graph\/})
of $E$, denoted
$\Ehat$ or $D(E)$, is the union of
$E$ and $E^*$. Thus, $\Ehat^0= E^0$ and $\Ehat^1= E^1\sqcup (E^1)^*$, with
range and source maps $r_{\Ehat}$ and $s_{\Ehat}$ combining those of $E$
and
$E^*$. When working with $\Ehat$, the edges (paths) from $E^1$ are often
called {\it real edges\/} ({\it real paths\/}) and those from $(E^1)^*$
{\it ghost edges\/} ({\it ghost paths\/}).
\enddefinition

\definition{1.3\. Path Algebras} Let $E$ be a graph and $K$ a field. The
{\it path algebra of $E$ over $K$\/}, denoted $KE$, is the $K$-algebra
based on the vector space over $K$ with basis the set of all paths in
$E$, and with multiplication induced from concatenation of paths: if $p=
e_1e_2\cdots e_n$ and $q= f_1f_2\cdots f_m$ are paths in $E$, their
product in $KE$ is given by
$$pq= \cases e_1e_2\cdots e_nf_1f_2\cdots f_m &\quad (\text{if\ } r(e_n)=
s(f_1)) \\  0 &\quad \text{(otherwise)}. \endcases$$

We shall need the observation that $KE$ is the $K$-algebra presented by
generators from the set $E^0\sqcup E^1$ with the following relations:
\roster
\item"\bull" $v^2=v$ for all $v\in E^0$;
\item"\bull" $vw=0$ for all distinct $v,w\in E^0$;
\item"\bull" $e=s(e)e= er(e)$ for all $e\in E^1$.
\endroster

The path algebra $KE$ is unital if and only if $E^0$ is finite, in which
case the identity is the sum of the vertices (= paths of length zero) in
$E$. In general, $KE$ is a {\it ring with local units\/}, meaning that
$KE$ contains a set $Q$ of pairwise commuting idempotents such that for
each $x\in KE$, there is some $q\in Q$ with $qx=xq=x$. Namely, take $Q$
to be the set of all finite sums of distinct vertices of $E$.
\enddefinition

\definition{1.4\. Leavitt Path Algebras} Before defining these algebras,
we recall two standard graph-theoretic concepts. A vertex $v$ in a graph
$E$ is a {\it sink\/} if $v$ emits no edges, i.e., there are no edges
$e\in E^1$ with $s(e)=v$. At the other extreme, $v$ is an {\it infinite
emitter\/} if there are infinitely many edges $e\in E^1$ with $s(e)=v$.

The {\it Leavitt path algebra\/} of $E$ over a field $K$,
denoted
$L_K(E)$, is the quotient of the path algebra $K\Ehat$ modulo the ideal
generated by the following elements:
\roster
\item"\bull" $e^*e-r(e)$, for all $e\in E^1$;
\item"\bull" $e^*f$, for all distinct $e,f\in E^1$;
\item"\bull" $v-\sum_{e\in E^1,\, s(e)=v} ee^*$, for every $v\in E^0$
which is neither a sink nor an infinite emitter.
\endroster
The corresponding equations in $L_K(E)$ are known as the {\it
Cuntz-Krieger relations\/}. Note that $L_K(E)$ is a ring with local units.

It is standard practice to use the same names for vertices, edges, and
paths in $\Ehat$ as for their cosets in $L_K(E)$. This introduces no
ambiguities when working with real paths or ghost paths (see Lemmas 1.5
and 1.6 below), but care must be taken with paths involving both real and
ghost edges. For instance, if $e\in E^1$, then $e^*e$ denotes a path of
length 2 in $K\Ehat$, while $e^*e= r(e)$ in $L_K(E)$. Let us reserve the
symbol
$\pi_{K,E}$ for the quotient map $K\Ehat\rightarrow L_K(E)$, for
use when it is important to distinguish cosets from their representatives.

Since the Cuntz-Krieger relations can be used to reduce any expression
involving a product of a ghost path followed by a real path, all elements
of $L_K(E)$ can be written as $K$-linear combinations of products $pq^*$
where $p$ and $q$ are real paths in $E$.

The path algebra $K\Ehat$ supports a $\ZZ$-grading under which real edges
have degree $1$ while ghost edges have degree $-1$. Thus, a real path $p$
is homogeneous of degree $\length(p)$, while a ghost path $p^*$ is
homogeneous of degree $-\length(p)$. The relations used to form $L_K(E)$
from $K\Ehat$ are all homogeneous (of degree zero), and therefore $L_K(E)$
inherits an induced $\ZZ$-grading.
\enddefinition

It has been noted in many papers that the vertices of a graph
are linearly independent when viewed as elements in a Leavitt path
algebra, but this fact does not appear to have been explicitly proved in
the literature. We take the opportunity to do so here.

\proclaim{Lemma 1.5} Let $K$ be a field and $E$ a graph. The
quotient map $\pi= \pi_{K,E} : K\Ehat \rightarrow L_K(E)$ sends the
vertices from $E^0$ to $K$-linearly independent elements
of $L_K(E)$. \endproclaim

\demo{Proof} We proceed by building a representation of
$L_K(E)$ as linear transformations on a vector space.

Let $\aleph$ be an infinite cardinal at least as large as
$\card(E^0\sqcup E^1)$. Let $X$ be an $\aleph$-dimensional vector space
over $K$, and set $R= \End_K(X)$. Since $\aleph\cdot \card(E^0) =
\aleph$, we can choose a decomposition $X= \bigoplus_{v\in E^0} X_v$ with
$\dim X_v =\aleph$ for all $v$. For $v\in E^0$, let $p_v\in R$ denote the
projection of $X$ onto $X_v$ with kernel $\bigoplus_{w\ne v} X_w$. For
each $v\in E^0$ which is not a sink, we can choose a decomposition $X_v=
\bigoplus_{e\in E^1,\, s(e)=v} Y_e$ with $\dim Y_e= \aleph$ for all $e$.
For $e\in E^1$, let $q_e\in R$ denote the projection of $X$ onto $Y_e$
with kernel 
$$(1-p_{s(e)})X \oplus \bigoplus \Sb f\in E^1,\, f\ne e\\ s(f)=s(e)
\endSb Y_f.$$ 
Finally, for each $e\in E^1$, choose $\alpha_e \in
q_eRp_{r(e)}$ and
$\alpha^*_e \in p_{r(e)}Rq_e$ such that $\alpha_e$ restricts to an
isomorphism $X_{r(e)} \rightarrow Y_e$ and $\alpha^*_e$ restricts to the
inverse isomorphism.

Observe that $p_v^2= p_v$ and $p_vp_w=0$ for all distinct $v,w\in E^0$.
For all $e\in E^1$, we have $\alpha_e= p_{s(e)}\alpha_e=
\alpha_ep_{r(e)}$ and $\alpha^*_e= p_{s(e^*)}\alpha^*_e=
\alpha^*_ep_{r(e^*)}$. Consequently, there is a unique $K$-algebra
homomorphism $\phi: K\Ehat \rightarrow R$ such that
\roster
\item"\bull" $\phi(v)= p_v$ for all $v\in E^0$;
\item"\bull" $\phi(e)= \alpha_e$ and $\phi(e^*)= \alpha^*_e$ for all $e\in
E^1$.
\endroster

We next check that $\ker\phi$ contains the defining relations of
$L_K(E)$. First, given $e\in E^1$, observe that $\alpha^*_e\alpha_e=
p_{r(e)}$, whence $e^*e-r(e) \in \ker\phi$. If $e,f\in E^1$ are distinct,
then $\alpha^*_e\alpha_f= \alpha^*_eq_eq_f\alpha_f= 0$, whence $e^*f\in
\ker\phi$. Finally, if $v\in E^0$ is neither a sink nor an infinite
emitter, then
$$p_v= \sum_{e\in E^1,\, s(e)=v} q_e= \sum_{e\in E^1,\, s(e)=r}
\alpha_e\alpha^*_e \,,$$
whence $v- \sum_{e\in E^1,\, s(e)=r} ee^* \in \ker\phi$. Thus, $\phi$
induces a unique $K$-algebra homomorphism $\psi: L_K(E) \rightarrow R$
such that $\psi\pi = \phi$. 

By construction, the projections $p_v$ for $v\in
E^0$ are pairwise orthogonal nonzero idempotents, and hence $K$-linearly
independent elements of
$R$. Since
$\psi\pi(v)= p_v$ for all $v$, we conclude that $\pi$ indeed maps the
elements of
$E^0$ to $K$-linearly independent elements of $L_K(E)$.  \qed\enddemo

That the result of Lemma 1.5 extends to (real) paths is proved in
\cite{\Sil, Lemma 1.1}. It is easy to include ghost paths, as follows.

\proclaim{Lemma 1.6} Let $K$ be a field and $E$ a graph. The quotient map
$\pi= \pi_{K,E} : K\Ehat \rightarrow L_K(E)$ restricts to an embedding of
the subspace $KE+KE^*$ of $K\Ehat$ into $L_K(E)$.
\endproclaim

\demo{Proof} For purposes of this proof, let us write $\xbar= \pi(x)$ for
$x\in K\Ehat$. We need to show that if 
$$\sum_{i=1}^m \alpha_i \pibar+ \sum_{j=1}^n \beta_j \qjstarbar =0
\tag\dagger$$  
for distinct real paths
$p_i$, distinct ghost paths $q^*_j$ of positive length, and scalars
$\alpha_i,\beta_j\in K$, then $\alpha_i=\beta_j=0$ for all $i$, $j$. If
$n=0$ and all the $p_i$ have length zero, this follows from Lemma 1.5.

Since $L_K(E)$ is a $\ZZ$-graded algebra, the equation $(\dagger)$ breaks
into homogeneous components, which can be treated separately. Hence, we
may assume that either $m=0$ or $n=0$, that all the $p_i$ have the same
length, and that all the $q_j$ have the same length. Further, $(\dagger)$
breaks into separate equations of similar form when the terms are
multiplied on the left by $\vbar$ and on the right by
$\wbar$, for any vertices $v,w\in E^0$. Consequently, it suffices to
consider the case where all the $p_i$ and $q^*_j$ have the same source
$v$ and the same range $w$.

Next, assume that $n=0$. The short argument of \cite{\Sil, Lemma 1.1}
applies here; we repeat it for the reader's convenience. Since the $p_i$
are distinct paths of the same length, from $v$ to $w$, we have
$\pjstarbar\pibar= \delta_{ij} \wbar$ for all $i$, $j$. Hence,
$\alpha_i \wbar =0$ for all $i$. However, $\wbar\ne 0$ by Lemma 1.5, and
therefore all $\alpha_i=0$.

The case $m=0$ is handled in the same manner. Namely, since the $q_j$ are
distinct paths of the same length, from $w$ to $v$, we have
$\qjstarbar\qbar_i= \delta_{ij} \vbar$ for all $i$, $j$. It follows that
all $\beta_j=0$, completing the proof. \qed\enddemo

\definition{1.7\. Opposites} The opposite algebras of path algebras and
Leavitt path algebras are easily understood in terms of dual graphs. If
$K$ is a field, $E$ a graph, and the notation $*$ is extended from edges
to paths as in (1.2), then $*$ provides a bijection from the standard
basis of $KE$ (the set of paths in $E$) onto the standard basis of $KE^*$
(the set of paths in $E^*$). Consequently, $*$ extends uniquely to a
$K$-vector space isomorphism of $KE$ onto $KE^*$. Since $*$ reverses the
composition of edges in paths, it also reverses multiplication of paths.
Therefore $*: KE\rightarrow KE^*$ is a $K$-algebra anti-isomorphism. In
particular, it follows that $KE^*$ is isomorphic to the opposite algebra
of $KE$.

The graph $\Ehat$ is self-dual, in that $\Ehat \cong \Ehat^*$, by an
isomorphism that fixes vertices and sends $e\mapsto (e^*)^*$ and
$e^*\mapsto e^*$ for $e\in E^1$. It is convenient to treat this
isomorphism as an identification, by setting
$(e^*)^* =e$ for all real edges $e$. Extending this to paths in the
natural manner, we obtain a $K$-algebra
anti-automorphism of
$K\Ehat$, still denoted $*$. In particular, $K\Ehat$ is isomorphic to its
opposite algebra.

Since the set of generators for the kernel of the quotient map $\pi_{K,E}$
given in (1.4) is mapped onto itself by $*$, this ideal is invariant under
$*$. Therefore
$*$ induces a $K$-algebra
anti-automorphism of $L_K(E)$, which we also denote $*$. Hence, $L_K(E)$
is isomorphic to its opposite.
\enddefinition

\definition{1.8\. Reversing the Path Composition Convention} At the level
of path algebras, reversing the convention for writing and composing
paths just changes the algebra to its opposite, but this is not so for
Leavitt path algebras. Just for the present subsection, let us
write $K{\updownarrow}E$ to denote the path algebra of $E$ over $K$
constructed using the right to left convention for paths. Thus, a path in
$K{\updownarrow}E$ is a product $p=e_ne_{n-1}\cdots e_1$ where
$e_1,\dots,e_n$ are edges from $E^1$ such that $r(e_i)= s(e_{i+1})$ for
$i=1,\dots,n-1$. Such a path runs from $s(p)=s(e_1)$ to $r(p)=r(e_n)$.
Multiplication of paths $p$ and
$q$ in $K{\updownarrow} E$ follows the rule that $pq=0$ unless $r(q)=s(p)$, in
which case $pq$ is the concatenation of $p$ and $q$ in the order ``first
$q$, then $p$''. Observe that $K{\updownarrow} E$ is the $K$-algebra presented by
generators from the set $E^0\sqcup E^1$ with the following relations:
\roster
\item"\bull" $v^2=v$ for all $v\in E^0$;
\item"\bull" $vw=0$ for all distinct $v,w\in E^0$;
\item"\bull" $e=r(e)e= es(e)$ for all $e\in E^1$.
\endroster
It follows that $K{\updownarrow} E$ equals the opposite
algebra of $KE$, and so $K{\updownarrow} E \cong KE^*$.

The construction of Leavitt path algebras requires choices of two
conventions -- one for paths, and one for the Cuntz-Krieger relations. If
we change both conventions from the ones in (1.1) and (1.4), we obtain
the quotient of $K{\updownarrow} E$ modulo the ideal generated by
\roster
\item"\bull" $ee^*-r(e)$, for all $e\in E^1$;
\item"\bull" $fe^*$, for all distinct $e,f\in E^1$;
\item"\bull" $v-\sum_{e\in E^1,\, s(e)=v} e^*e$, for every $v\in E^0$
which is neither a sink nor an infinite emitter.
\endroster
This algebra is just the opposite algebra of $L_K(E)$, and so it is
isomorphic to $L_K(E)$.

However, if we adopt the right to left convention for paths without
reversing products in the Cuntz-Krieger relations, the latter, to be
sensible, must be modified by interchanging sources and ranges. This time,
we obtain the quotient of
$K{\updownarrow} E$ modulo the ideal generated by
\roster
\item"\bull" $e^*e-s(e)$, for all $e\in E^1$;
\item"\bull" $e^*f$, for all distinct $e,f\in E^1$;
\item"\bull" $v-\sum_{e\in E^1,\, r(e)=v} ee^*$, for every $v\in E^0$
which is neither a source nor an infinite receiver.
\endroster
Since $s(e)= r(e^*)$ and $r(e)= s(e^*)$ for $e\in E^1$, the opposite of
the above algebra is naturally isomorphic to $L_K(E^*)$, and hence the
algebra itself is isomorphic to $L_K(E^*)$. Thus, to transfer results
proved about Leavitt path algebras constructed with these conventions, one
must dualize any graph conditions that appear.
\enddefinition

Much of the initial work on Leavitt path algebras restricted attention to
graphs of the following type.

\definition{1.9\. Row-Finite Graphs} The {\it incidence matrix\/} for a
graph $E$ is an $E^0\times E^0$ matrix in which the $(v,w)$-entry (for
vertices $v,w\in E^0$) is the number of edges from $v$ to $w$ in
$E^1$. We say that $E$ is a {\it row-finite graph\/} provided its
incidence matrix is row-finite with finite entries, i.e., each row
contains at most finitely many nonzero entries, none of which is $\infty$.
Thus, $E$ is row-finite if and only if each vertex of $E$ emits at most
finitely many edges.
\enddefinition

\head 2. Direct Limits \endhead

Here we set up an appropriate category of graphs on which the
construction of Leavitt path algebras (over a given field) is functorial,
and show that this functor preserves direct limits. It is helpful to base
the discussion on direct limits of ordinary path algebras, and so we
begin with functoriality of the path algebra construction. Recall that the
term {\it direct limit\/} (equivalently, {\it inductive limit\/}) refers
to a colimit (in a category) over a system of objects and morphisms
indexed by a {\it directed\/} set.

\definition{2.1\. The Directed Graph Category} A {\it graph morphism\/}
from a graph $E$ to a graph $F$ is a pair $\phi= (\phi^0,\phi^1)$
consisting of maps $\phi^0: E^0\rightarrow F^0$ and
$\phi^1: E^1\rightarrow F^1$ such that $s_F\phi^1= \phi^0s_E$ and
$r_F\phi^1= \phi^0r_E$. Due to the assumption that $E^0\cap E^1=F^0\cap
F^1= \varnothing$, we can view $\phi$ as a function $E^0\sqcup
E^1\rightarrow F^0\sqcup F^1$ that restricts to $\phi^0$ and $\phi^1$.

Particularly useful graph morphisms arise when $E$ is a {\it subgraph\/}
of $F$, meaning that $E$ consists of some of the vertices and edges of
$F$. Of course, a collection of vertices and edges from $F$ does not
naturally form a graph unless the source and range vertices of the chosen
edges are included among the chosen vertices. Thus, to say that $E$ is a
subgraph of $F$ means that
\roster
\item"\bull" $E^0 \subseteq F^0$ and $E^1 \subseteq F^1$;
\item"\bull" $r_F(e),s_F(e)\in E^0$ for all $e\in E^1$;
\item"\bull" $r_E= r_F|_{E^1}$ and $s_E= s_F|_{E^1}$.
\endroster
When $E$ is a subgraph of $F$, the inclusion map $E^0\sqcup
E^1\rightarrow F^0\sqcup F^1$
forms a graph morphism
$E\rightarrow F$.

We shall let $\DiGr$ denote the {\it category of directed graphs\/}: the
objects of $\DiGr$ are arbitrary directed graphs, and the morphisms are
arbitrary graph morphisms.
\enddefinition

The construction of path algebras over a field $K$ appears functorial at
first glance, but there is one problem with preservation of relations.
For example, consider the graph morphism
$$\phi: E = \bigl( {\xymatrix{ v \ar[r]^e &w }} \bigr) \longrightarrow F=
\biggl( {\xymatrix{ x \ar@(dr,ur)_(0.75)f }}\ \biggr)$$ 
sending the vertices $v,w\in E^0$ to the vertex $x\in F^0$, and the
edge $e\in E^1$ to the edge $f\in F^1$. Since $vw=0$ in $KE$ while
$x^2=x\ne 0$ in
$KF$, there is no
$K$-algebra morphism $KE \rightarrow KF$ extending $\phi$. To avoid this
problem, we restrict attention to graph morphisms which are injective on
vertices.

\definition{2.2\. Path Algebra Functors} Let $\DiGr_0$ be the
subcategory of $\DiGr$ whose objects are arbitrary directed graphs and
whose morphisms are those graph morphisms $\phi$ for which $\phi^0$ is
injective. Given a field $K$, let $\KAlg$
denote the category of (not necessarily unital) $K$-algebras. We allow
completely arbitrary $K$-algebras as objects of $\KAlg$ (that is,
arbitrary vector spaces over $K$, equipped with associative, $K$-bilinear
multiplications), and arbitrary $K$-algebra morphisms between them (that
is, arbitrary multiplicative $K$-linear maps). Any direct system in
$\KAlg$ has a direct limit in this category.

Now if $\phi$ is a morphism in $\DiGr_0$, we have
\roster
\item"\bull" $\phi^0(v)^2=\phi^0(v)$ in $KF$ for all $v\in E^0$
(because
$\phi^0$ sends vertices to vertices);
\item"\bull" $\phi^0(v)\phi^0(w)=0$ in $KF$ for all distinct $v,w\in
E^0$ (because $\phi^0$ sends distinct vertices to distinct vertices);
\item"\bull" $\phi^1(e)= \phi^0s(e)\phi^1(e)= \phi^1(e)\phi^0r(e)$ in
$KF$ for all $e\in E^1$ (because $\phi^0s(e)= s(\phi^1(e))$ and
$\phi^0r(e)= r(\phi^1(e))$).
\endroster
Consequently, the map $\phi: E^0\sqcup E^1 \rightarrow
F^0\sqcup F^1$ uniquely extends to a $K$-algebra morphism $K\phi:
KE\rightarrow KF$.

The assignments $E\mapsto KE$ and $\phi\mapsto K\phi$ define a functor
$K[{-}]: \DiGr_0 \rightarrow \KAlg$.
\enddefinition

There are additional problems with the construction of Leavitt
path algebras over a field $K$, because the induced morphisms between
path algebras of doubles do not always preserve Cuntz-Krieger relations.
For instance, consider the graph morphism
$$\psi: E= \bigl( {\xymatrix{ v \ar[r]^e &w }} \bigr)
\longrightarrow F= \bigl( {\xymatrix{ v' \ar@/^/[r]^f \ar@/_/[r]_g &w' }}
\bigr)$$ 
sending $v\mapsto v'$ and $w\mapsto w'$, while $e\mapsto f$. Now $ee^*=v$
in $L_K(E)$, while $ff^*= v'-gg^*\ne v'$ in $L_K(F)$ (e.g., $g^*(gg^*)g=
w'\ne 0$, so $gg^*\ne 0$). Thus, there is no
$K$-algebra morphism $L_K(E) \rightarrow L_K(F)$ extending $\psi$.
Similarly, the graph morphism
$$F= \bigl( {\xymatrix{ v' \ar@/^/[r]^f \ar@/_/[r]_g &w' }}
\bigr) \longrightarrow E= \bigl( {\xymatrix{ v \ar[r]^e &w }} \bigr)$$
sending $v'\mapsto v$ and $w'\mapsto w$, while $f,g\mapsto e$, does not
extend to a $K$-algebra morphism $L_K(F) \rightarrow L_K(E)$, because
$f^*g=0$ in $L_K(F)$ while $e^*e= w\ne 0$ in $L_K(E)$.

These problems can be dealt with by further restricting the allowable
morphisms between graphs.

\definition{2.3\. The Cuntz-Krieger Graph Category} Let us say that a
graph morphism $\phi: E\rightarrow F$ is a {\it CK-morphism\/} (short for
{\it Cuntz-Krieger morphism\/}) provided
\roster
\item $\phi^0$ and $\phi^1$ are both injective;
\item For each $v\in E^0$ which is neither a sink nor an infinite emitter,
$\phi^1$ induces a bijection $s_E^{-1}(v) \rightarrow
s_F^{-1}(\phi^0(v))$. 
\endroster
In particular, condition (1) says
that $\phi$ maps $E$ isomorphically onto a subgraph of $F$, while
condition (2) implies that
$\phi^0$ must send non-sink finite emitters to non-sink finite emitters.

If $E$ is row-finite, injectivity of $\phi^0$ together with condition (2)
is sufficient to ensure injectivity of $\phi^1$. Thus, in this case,
$\phi$ is a CK-morphism if and only if it is a {\it complete graph
homomorphism\/} in the sense of \cite{\AMP, p\. 161}.

We shall say that a subgraph $E$ of a graph $F$ is a {\it CK-subgraph\/}
provided the inclusion map $E\rightarrow F$ is a CK-morphism.

Now define $\CKGr$ to be the subcategory of $\DiGr$ whose objects are
arbitrary directed graphs and whose morphisms are arbitrary CK-morphisms.
\enddefinition

\definition{2.4\. Leavitt Path Algebra Functors} Fix a field $K$, and
consider a CK-morphism $\phi: E\rightarrow F$ between graphs $E$ and $F$.
Now $\phi$ extends to a graph morphism $\phihat: \Ehat\rightarrow \Fhat$
sending $e^*\mapsto \phi^1(e)^*$ for all $e\in E^1$, and $\phihat^0=
\phi^0$ is injective by assumption. Hence, $\phihat$ uniquely induces a
$K$-algebra morphism $K\phihat: K\Ehat\rightarrow K\Fhat$. Observe that
\roster
\item"\bull" $\phihat^1(e^*)\phihat^1(e)=
r(\phihat^1(e))= \phihat^0(r(e))$ in $L_K(F)$ for all $e\in E^1$;
\item"\bull" $\phihat^1(e^*)\phihat^1(f)= 0$ in
$L_K(F)$ for all distinct $e,f\in E^1$ (because $\phihat^1(e)$ and
$\phihat^1(f)$ are distinct edges in $F$);
\item"\bull" For every $v\in E^0$ which is neither a sink nor an
infinite emitter, the same properties hold for $\phihat^0(v)$, and
$\phihat^0(v)= \sum_{f\in F^1,\, s(f)=\phihat^0(v)} ff^*=
\sum_{e\in E^1,\, s(e)=v} \phihat^1(e)\phihat^1(e^*)$ in $L_K(F)$ (because
of condition (2.3)(2)).
\endroster
Consequently, $K\phihat$ uniquely induces a $K$-algebra morphism
$L_K(\phi): L_K(E) \rightarrow L_K(F)$.

The assignments $E\mapsto L_K(E)$ and $\phi\mapsto L_K(\phi)$ define a
functor $L_K: \CKGr\rightarrow \KAlg$.
\enddefinition

\proclaim{Lemma 2.5} {\rm (a)} Arbitrary direct limits exist in the
categories $\DiGr$, $\DiGr_0$, and $\CKGr$.

{\rm (b)} For any field $K$, the functors $K[{-}]: \DiGr_0
\rightarrow \KAlg$ and $L_K: \CKGr\rightarrow
\KAlg$ preserve direct limits.
\endproclaim

\demo{Proof} (a) Let $\E= \bigl( (E_i)_{i\in I},(\phi_{ij})_{i\le j\,
\text{in}\, I} \bigr)$ be a direct system in $\DiGr$. Let
$E^0_\infty$ and $E^1_\infty$ denote the direct limits of the
corresponding direct systems of sets, $\bigl( (E^0_i)_{i\in
I},(\phi^0_{ij})_{i\le j\, \text{in}\, I} \bigr)$ and $\bigl(
(E^1_i)_{i\in I},(\phi^1_{ij})_{i\le j\, \text{in}\, I} \bigr)$, with
limit maps
$\eta^0_i: E^0_i \rightarrow E^0_\infty$ and $\eta^1_i: E^1_i \rightarrow
E^1_\infty$. 

Given $f\in E^1_\infty$, the set $I_f= \{i\in I\mid f\in
\eta^1_i(E^1_i)\}$ is nonempty and upward directed. For any $i,j\in I_f$
and any $e\in (\eta^1_i)^{-1}(\{f\})$ and
$e'\in (\eta^1_j)^{-1}(\{f\})$, we have $\eta^1_i(e)=f= \eta^1_j(e')$,
and so there is an index $k\ge i,j$ in $I$ such that $\phi^1_{ik}(e)=
\phi^1_{jk}(e')$. Consequently,
$$\eta^0_i(r(e))= \eta^0_k\phi^0_{ik}(r(e))= \eta^0_k\bigl(
r(\phi^1_{ik}(e)) \bigr)= \eta^0_k\bigl(
r(\phi^1_{jk}(e')) \bigr)= \eta^0_j(r(e')).$$
Thus, there is a unique well-defined map $r_\infty: E^1_\infty \rightarrow
E^0_\infty$ such that $r_\infty(\eta^1_i(e))= \eta^0_i(r(e))$ for all
$i\in I$ and $e\in E^1_i$. Similarly, there is a unique well-defined map
$s_\infty: E^1_\infty \rightarrow E^0_\infty$ such that
$s_\infty(\eta^1_i(e))= \eta^0_i(s(e))$ for all
$i\in I$ and $e\in E^1_i$.

Now $E_\infty= (E^0_\infty,E^1_\infty,r_\infty,s_\infty)$ is a directed
graph, and the maps $\eta_i= \eta^0_i\sqcup \eta^1_i: E_i \rightarrow
E_\infty$ are morphisms in $\DiGr$ such that $\eta_j\phi_{ij}= \eta_i$
for all $i\le j$ in $I$. It is routine to check that $E_\infty$, together
with the maps $\eta_i$, is a direct limit for the system $\E$ in $\DiGr$.

If all the maps $\phi^0_{ij}$ are injective, then so are the maps
$\eta^0_i$, and $E_\infty$ is a direct limit for $\E$ in $\DiGr_0$.

Now assume that all the $\phi_{ij}$ are CK-morphisms, and note that the
maps $\eta^0_i$ and $\eta^1_i$ are injective. Consider a vertex $v\in
E^0_\infty$ which is neither a sink nor an infinite emitter. Then
$v=\eta^0_i(w)$ for some
$i\in I$ and some $w\in E^0_i$. Since $\eta^1_i$ is injective, we see
that $w$ cannot be an infinite emitter. There is at least one edge $e\in
E^1_\infty$ emitted by $v$, and, after possibly increasing
$i$, we may assume that
$e=\eta^1_i(f)$ for some edge $f\in E^1_i$ emitted by $w$. Hence, $w$ is
not a sink. For all $j\ge i$ in $I$, the map $\phi^1_{ij}$ sends the
edges emitted by $w$ bijectively onto the edges emitted by
$\phi^0_{ij}(w)$, and consequently $\eta^1_i$ sends the edges emitted by
$w$ bijectively onto the edges emitted by $v$. Thus, the $\eta_i$ are
CK-morphisms. 

Suppose we are given a graph $F$ and CK-morphisms $\theta_i:
E_i\rightarrow F$ such that $\theta_j\phi_{ij}= \theta_i$ for all $i\le
j$ in $I$. Since $E_\infty$ is a direct limit for $\E$ in $\DiGr_0$, there
is at least a unique graph morphism
$\sigma: E_\infty
\rightarrow F$ such that $\sigma\eta_i= \theta_i$ for all $i\in I$. In
the same manner as above, one checks that $\sigma$ is a CK-morphism.
Therefore $E_\infty$ is a direct limit for $\E$ in $\CKGr$.

(b) Let $\E= \bigl( (E_i)_{i\in I},(\phi_{ij})_{i\le j\,
\text{in}\, I} \bigr)$ be a direct system in $\DiGr_0$, and construct
the graph  $E_\infty$ and maps
$\eta_i$ as in part (a). We first show that $KE_\infty$ is a direct limit
for the system $K[\E]= \bigl( (KE_i)_{i\in I},(K\phi_{ij})_{i\le j\,
\text{in}\, I} \bigr)$ in $\KAlg$.

Let a $K$-algebra $A$, with limit  maps $\theta_i: KE_i \rightarrow A$,
be a direct limit for the system $K[\E]$. Since $(K\eta_j)(K\phi_{ij})=
K\eta_i$ for all $i\le j$ in $I$, there is a unique $K$-algebra morphism
$\psi: A\rightarrow KE_\infty$ such that $\psi\theta_i= K\eta_i$ for all
$i\in I$. We show that $\psi$ is an isomorphism.

Any element $x\in KE_\infty$ is a $K$-linear combination of paths
involving finitely many edges from $E^1_\infty$. The edges needed can all
be found in 
 $\eta^1_i(E^1_i)$ for some $i\in
I$, whence $x\in (K\eta_i)(KE_i) \subseteq \psi(A)$. Consequently, $\psi$
is surjective. Now consider $a\in \ker\psi$, write $a=\theta_i(b)$ for 
 some $i\in I$ and $b\in KE_i$, and write $b= \sum_{l\in L} \alpha_l p_l$
for some finite index set $L$, some
$\alpha_l\in K$, and some paths $p_l\in E^1_i$. Then
$$\sum_{l\in L} \alpha_l (K\eta_i)(p_l)= \psi\theta_i(b) =0,$$
where the $(K\eta_i)(p_l)$ can be viewed as paths in $E_\infty$. It
follows that there is a partition $L= L_1\sqcup \cdots\sqcup L_t$ such
that for each
$m=1,\dots,t$, the paths $(K\eta_i)(p_l)$ for $l\in L_m$ are all equal,
and $\sum_{l\in L_m} \alpha_l =0$. There is an index $j\ge i$ in $I$ such
that for each
$m=1,\dots,t$, the paths $(K\phi_{ij})(p_l)$ for $l\in L_m$ are all
equal.  Hence, $(K\phi_{ij})(b)= \sum_{l\in L} \alpha_l (K\phi_{ij})(p_l)
=0$ and so $a= \theta_j(K\phi_{ij})(b) =0$. Thus $\psi$ is an
isomorphism, as announced.

Therefore, we have shown that $K[{-}]$ preserves direct limits.

Now assume that $\E$ is a direct system in $\CKGr$, that is, all the maps
$\phi_{ij}$ are CK-morphisms. As  noted above, the limit maps $\eta_i$
are now CK-morphisms. Taking doubles (which is functorial), we obtain a
direct system
$\calEhat= \bigl( (\Ehat_i)_{i\in I},(\phihat_{ij})_{i\le j\,
\text{in}\, I} \bigr)$, which lives in $\DiGr_0$. By what we have
proved so far, $K\Ehat_\infty$ is a direct limit for $K[\calEhat]$, with
limit maps $K\etahat_i : K\Ehat_i \rightarrow K\Ehat_\infty$. 

For $i\le j$ in $I$, we have $\pi_{K,E_j}(K\phihat_{ij})=
L_K(\phi_{ij})\pi_{K,E_i}$ and so
$(K\phihat_{ij})(\ker\pi_{K,E_i})
\subseteq \ker\pi_{K,E_j}$. Similarly,
$(K\etahat_i)(\ker\pi_{K,E_i}) \subseteq \ker\pi_{K,E_\infty}$ for all
$i\in I$. We claim that any of the standard generators of
$\ker\pi_{K,E_\infty}$ must lie in $(K\etahat_i)(\ker\pi_{K,E_i})$ for
some $i$. First, for any
$e\in E^1_\infty$, we have $e=\eta^1_i(e_i)$ for some $i\in I$ and some
$e_i\in E^1_i$, whence $(K\etahat_i)(e_ie^*_i-r(e_i))= ee^*-r(e)$.
Similarly, any distinct edges $e,f\in E^1_\infty$ arise as
$e=\eta^1_i(e_i)$ and $f=\eta^1_i(f_i)$ for some $i\in I$ and some
distinct edges $e_i,f_i\in E^1_i$, whence $(K\etahat_i)(e_if^*_i)=
ef^*$. Finally, consider a vertex $v\in E^0_\infty$ which is
neither a sink nor an infinite emitter. As in the proof of part (a),
there exist $i\in I$ and a vertex $w\in E^0_i$ such that $\eta^0_i(w)=v$
and $\eta^1_i$ sends the edges emitted by
$w$ bijectively onto the edges emitted by $v$. Thus,
$$v-\sum_{e\in
E^1_\infty,\, s(e)=v} ee^*= (K\etahat_i) \biggl( w-\sum_{f\in
E^1_i,\, s(f)=w} ff^* \biggr),$$
and the claim is proved. It follows that
$$\ker\pi_{K,E_\infty}= \bigcup_{i\in I}\,
(K\etahat_i)(\ker\pi_{K,E_i}).$$

Therefore $K\Ehat_\infty/\ker\pi_{K,E_\infty}$ is the direct limit of the
quotients $K\Ehat_i/\ker\pi_{K,E_i}$, that is, $L_K(E_\infty)$ is the
direct limit of the system $L_K(\E)$ in $\CKGr$. \qed\enddemo

\proclaim{Proposition 2.6} {\rm [Ara-Moreno-Pardo]} If $E$ is a row-finite
graph, then
$E$ is the direct limit {\rm(}in $\CKGr${\rm)} of its finite CK-subgraphs.
Consequently, the Leavitt path algebra of $E$, over any field $K$, is the
direct limit of the Leavitt path algebras of the finite CK-subgraphs of
$E$.  \endproclaim

\demo{Proof} As noted in (2.3), the CK-morphisms between row-finite
graphs coincide with the complete graph homomorphisms in the sense of
\cite{\AMP}. In particular, the CK-subgraphs of $E$ are just the {\it
complete subgraphs\/} of \cite{\AMP}. Hence, the proposition is just a
restatement of
\cite{\AMP, Lemmas 3.1, 3.2}.  \qed\enddemo 

In general, Proposition 2.6 fails for row-infinite graphs, since such a
graph may have too meager a supply of finite CK-subgraphs. For instance,
take the graph $E$ consisting of one vertex $v$ and infinitely many edges
(all of which must be loops from
$v$ to $v$); the only finite CK-subgraphs of $E$ are the empty graph
(with no vertices or edges) and the graph with one vertex $v$ but no
edges. However, there is a general result for direct limits of countable
subgraphs, as follows.

\proclaim{Proposition 2.7} Let $E$ be an arbitrary graph, and $K$ a
field. Then $E$ is the direct limit of its countable CK-subgraphs, and
consequently $L_K(E)$ is the direct limit of the $L_K(F)$ over countable
CK-subgraphs $F$ of $E$. \endproclaim

\demo{Proof} For the first conclusion, we just need to show that $E$ is
the directed union of its countable CK-subgraphs, i.e., each vertex or
edge of $E$ lies in at least one countable CK-subgraph, and any two
countable CK-subgraphs of $E$ are contained in a common countable
CK-subgraph. Once this is established, the second conclusion follows via
Lemma 2.5(b). The required properties are both consequences of the
following claim:
\roster
\item"\bull" {\bf Claim.} Given any countable subsets $X^0\subseteq E^0$
and $X^1 \subseteq E^1$, there exists a countable CK-subgraph $F$ of $E$
such that $X^0\subseteq F^0$ and $X^1\subseteq F^1$.
\endroster

Starting with $X^0$ and $X^1$ as in the claim, we construct a countable
ascending sequence of countable subgraphs of $E$, labelled $F_0\subseteq
F_1\subseteq \cdots$. To begin, let $F_0$ be the subgraph of $E$
generated by $X^0\sqcup X^1$, that is, $F^1_0= X^1$ and
$F^0_0= X^0 \sqcup r_E(X^1) \cup s_E(X^1)$, with $r_{F_0}$ and $s_{F_0}$
being the restrictions of $r_E$ and $s_E$ to $F^1_0$.

Once $F_n$ has been constructed, choose sets $X^1_n(v) \subseteq E^1$ for
each vertex $v\in F^0_n$ as follows:
\roster
\item"\bull" If $v$ emits at most finitely many edges in $E$, set
$X^1_n(v)= s_E^{-1}(v)$;
\item"\bull" If $v$ emits infinitely many edges in $E$, take $X^1_n(v)$
to be some countably infinite subset of $s_E^{-1}(v)$.
\endroster
Then, set $F_{n+1}$ equal to the subgraph of $E$ generated by $F^0_n\sqcup
F^1_n\cup \bigcup_{v\in F^0_n} X^1_n(v)$.

Finally, $F= \bigcup_{n=0}^\infty F_n$ is a countable CK-subgraph of $E$
containing $X^0\sqcup X^1$. \qed\enddemo

\definition{2.8\. Some Notation} Since we will utilize the direct limits
given in Proposition 2.7 a number of times, it is convenient to establish
some corresponding notation. Given a graph $E$, write
$(E_\alpha)_{\alpha\in A}$ for the family of all countable CK-subgraphs
of $E$. Inclusions among these subgraphs translate into a partial
ordering on the index set $A$, where $\alpha,\beta \in A$ satisfy
$\alpha\le \beta$ if and only if $E_\alpha \subseteq E_\beta$. As shown
in the proof of Proposition 2.7, any countable subgraph of $E$ is
contained in some $E_\alpha$. In particular, this means that $A$ is
countably upward directed: given any countable sequence
$\alpha_1,\alpha_2,\dots$ in $A$, there is some $\alpha\in A$ such that
all $\alpha_i \le \alpha$. Let us also observe that any nonempty countable
upward directed subset $B\subseteq A$ has a supremum in $A$. To see this,
we just need to observe that the union of the subgraphs indexed by $B$ is
a countable CK-subgraph of $E$. In particular, any ascending sequence
$\alpha_1\le \alpha_2\le \cdots$ in $A$ has a supremum.

For $\alpha\in A$, let $\eta_\alpha$ denote the inclusion map $E_\alpha
\rightarrow E$, and for $\alpha\le \beta$ in $A$, let
$\phi_{\alpha\beta}$ denote the inclusion map $E_\alpha \rightarrow
E_\beta$. Then we have
$$E= \dirlim \bigl( (E_\alpha)_{\alpha\in A}\,,\,
(\phi_{\alpha\beta})_{\alpha\le\beta\,\text{in}\,A} \bigr),$$
with limit maps $\eta_\alpha$. For any field $K$, we then have
$$L_K(E)= \dirlim \bigl( \bigl( L_K(E_\alpha) \bigr)_{\alpha\in A}\,,\,
\bigl( L_K(\phi_{\alpha\beta}) \bigr)_{\alpha\le\beta\,\text{in}\,A}
\bigr),$$ 
with limit maps $L_K(\eta_\alpha)$.
\enddefinition

\head 3. Ideals\endhead

A number of key results about the ideal theory of Leavitt path algebras
of countable graphs are established in \cite{\Tom} and \cite{\AAc}.
While some of the proofs work equally well for uncountable graphs, the 
main results rely on a process called {\it desingularization\/},
introduced in \cite{\DT}, which only applies to graphs in which each
vertex emits at most countably many edges. A desingularization of a
countable graph $E$ is a countable row-finite graph $E'$, such that for
any field $K$, the algebras $L_K(E)$ and $L_K(E')$ are Morita equivalent
\cite{\AAc, Theorem 5.2} (see \cite{\Tom, Lemma 6.7} for another proof).
Since the proof of this theorem relies on the countability of $E^0$,
there is at present no usable theory of desingularization for uncountable
graphs. The general theorems, however, can be extended to uncountable
graphs by direct limit techniques, as we show below. We begin with a
rough outline of a key point of the procedure.

\definition{3.1\. Modus Operandi} Let $E$ be a graph, expressed as a
direct limit of its countable CK-subgraphs as in (2.8). For suitable
``finitely definable'' graph-theoretic properties $\calP$ (a term we do not
make completely precise), a sufficient supply of countable CK-subgraphs of
$E$ satisfying
$\calP$ can be constructed as follows.

Assume that $\calP$ can be expressed in the form ``for any choice of
finitely many vertices and edges satisfying a certain finite list of
conditions, there exist finitely many vertices and edges satisfying
another given finite list of conditions''. If $E$ satisfies $\calP$, then
given any finite subset $X$ of $E$ (that is, $X\subseteq E^0\sqcup E^1$)
satisfying the hypotheses of $\calP$, there is a finite subgraph $F$ of $E$
that contains
$X$ together with a finite subset satisfying the conclusions of $\calP$.
Given any index $\alpha\in A$, the countable subgraph $E_\alpha$ has
only countably many finite subsets, and so $E$ has a countable subgraph
$G$ that contains $E_\alpha$ along with finite subsets satisfying the
conclusions of $\calP$ relative to any finite subset of $E_\alpha$. Since
every countable subgraph of $E$ is contained in a countable CK-subgraph,
$G\subseteq E_\beta$ for some $\beta\ge\alpha$ in $A$. Now repeat this
procedure countably many times, obtaining indices $\beta(0)=\alpha \le
\beta(1)\le\cdots$ in $A$ such that for each finite subset $X$ of any
$E_{\beta(n)}$ satisfying the hypotheses of $\calP$, there is a finite
subset of $E_{\beta(n+1)}$ satisfying the conclusions of $\calP$ relative
to $X$. Finally, taking
$\beta$ to be the supremum of the $\beta(n)$ in $A$, we observe that
$E_\beta= \bigcup_{n=0}^\infty E_{\beta(n)}$ satisfies $\calP$.

Therefore, if the property $\calP$ behaves as required in the above outline,
we conclude that for each $\alpha\in A$, there is some $\beta\ge\alpha$
in $A$ such that $E_\beta$ satisfies $\calP$. Consequently, $E$ is a direct
limit in $\CKGr$ of countable graphs satisfying $\calP$.
\enddefinition

As a first illustration of our M.O., we consider the {\it homogeneous
ideals\/} of Leavitt path algebras $L_K(E)$, with respect to the
$\ZZ$-grading introduced in (1.4). (In many references, these ideals are
called {\it graded\/} rather than homogeneous.) Tomforde proved in
\cite{\Tom, Theorem 4.8} that if $E$ is countable, $\pi$ is a graded
ring homomorphism from $L_K(E)$ to a $\ZZ$-graded ring, and $\pi(v)\ne 0$
for all $v\in E^0$, then $\pi$ is injective. An equivalent conclusion is
that every nonzero homogeneous ideal of $L_K(E)$ has nonempty
intersection with $E^0$. Tomforde's proof, which does not use
desingularization, appears to work for uncountable graphs just as well as
for countable ones. That involves checking through a long sequence of
details, however. As direct limits offer a quick means to carry over the
result from the countable to the uncountable case, we prove it that way,
to show off the method.

\proclaim{Theorem 3.2}  {\rm (Extending \cite{\Tom, Theorem 4.8})} Let $K$
be a field and $E$ a graph. Then every nonzero homogeneous ideal of
$L_K(E)$ has nonempty intersection with $E^0$. \endproclaim

\demo{Proof} Write $E= \dirlim E_\alpha$ and $L_K(E)= \dirlim
L_K(E_\alpha)$ as in (2.8), and let $I$ be a nonzero homogeneous ideal of
$L_K(E)$. Each of the limit maps $L_K(\eta_\alpha)$ is a graded
homomorphism, whence $I_\alpha= L_K(\eta_\alpha)^{-1}(I)$ is a homogeneous
ideal of $L_K(E_\alpha)$. Further, $I=
\bigcup_{\alpha\in A}  L_K(\eta_\alpha)(I_\alpha)$, and so there must be
some $\beta\in A$ such that $I_\beta \ne 0$. By \cite{\Tom,
Theorem 4.8}, there exists a vertex $v\in I_\beta \cap E^0_\beta$,
whence $L_K(\eta_\beta)(v) \in I\cap E^0$. \qed\enddemo

\proclaim{Corollary 3.3} Let $K$ be a field. For any morphism
$\phi:E\rightarrow F$ in $\CKGr$, the $K$-algebra homomorphism $L_K(\phi):
L_K(E) \rightarrow L_K(F)$ is injective. \endproclaim

\demo{Proof} Since $L_K(\phi)$ is a graded ring homomorphism, its kernel
is a homogeneous ideal of $L_K(E)$. Moreover, $L_K(\phi)$ maps any vertex
$v\in E^0$ to a vertex $\phi^0(v)\in F^0$, and so $L_K(\phi)(v)\ne 0$ by
Lemma 1.5. Therefore Theorem 3.2 implies $\ker L_K(\phi) =0$. \qed\enddemo

In order to address arbitrary (non-homogeneous) ideals of Leavitt path
algebras, some graph-theoretic properties must be considered.

\definition{3.4\. Conditions (K) and (L)} Let $E$ be a graph. A path
$p=e_1e_2\cdots e_n$ in $E$ is {\it closed\/} if $r(p)=s(p)$, in which
case $p$ is said to be {\it based at\/} the vertex $r(p)=s(p)$. A closed
path $p$ as above is {\it simple\/} provided it does not pass through its
base more than once, i.e., $s(e_i)\ne s(e_1)=s(p)$ for all $i=2,\dots,n$.
More strictly, $p$ is a {\it cycle\/} if it is closed and it does not pass
through any vertex twice, i.e., $s(e_i)\ne s(e_j)$ for all $i\ne j$. An
{\it exit\/} for $p$ is an edge $e\in E^1$ that starts at some vertex on
$p$ but does not coincide with the next edge of $p$, i.e., there is an
index $i\in \{1,\dots,n\}$ such that $s(e)= s(e_i)$ but $e\ne e_i$.
(However, it is allowed that $r(e)$ might equal some $s(e_j)$.) Note that
if $p$ is a simple closed path which is not a cycle,
then $p$ automatically has an exit: There are indices $1<i<j\le n$ such
that $s(e_i)=s(e_j)$ but $e_i\ne e_j$, so $e_j$ is an exit at
$s(e_i)$. 

The graph $E$ satisfies {\it Condition\/} (K) provided no vertex $v\in
E^0$ is the base of precisely one simple closed path, i.e., either no
simple closed paths are based at $v$, or at least two are based there. It
satisfies {\it Condition\/} (L) provided every simple closed path in $E$
has an exit, or, equivalently, every cycle in $E$ has an exit.

Observe that (K)$\implies$(L). For if $p=e_1e_2\cdots e_n$ is a cycle in
$E$, then Condition (K) implies that $E$ contains a simple closed path
$q= f_1f_2\cdots f_m \ne p$ based at $s(p)$. We cannot have $e_i=f_i$ for
$i\le \min\{m,n\}$, because then the longer of $p$ or $q$ would pass
through $s(p)=s(q)$ more than once. Take $j\le \min\{m,n\}$ to be the
least index such that $e_j\ne f_j$. Then $s(e_j)=s(f_j)$ and $f_j$ is an
exit for $p$.
\enddefinition

\proclaim{Lemma 3.5} Let $E$ be a graph that satisfies Condition {\rm(K)
(}respectively, Condition {\rm(L))}. Write $E= \dirlim E_\alpha$ as in
{\rm(2.8)}. For each $\alpha\in A$, there exists $\beta\ge\alpha$ in $A$
such that $E_\beta$ satisfies Condition {\rm(K)
(}respectively, Condition {\rm(L))}. Consequently, $E$ is a direct limit
in $\CKGr$ of countable graphs satisfying Condition {\rm(K)
(}respectively, Condition {\rm(L))}. \endproclaim

\demo{Proof} Assume first that $E$ satisfies Condition (K). Given
$\alpha\in A$, let $v_1,v_2,\dots$ be a list of those vertices in
$E^0_\alpha$ which are bases for simple closed paths in $E_\alpha$. For
$i=1,2,\dots$, there is a simple closed path $p_i$ in $E_\alpha$ based at
$v_i$, and since $E$ satisfies Condition (K), there must be a simple
closed path $q_i\ne p_i$ in $E$ which is based at $v_i$. There is some
$\gamma\ge\alpha$ in $A$ such that $E_\gamma$ contains all the $q_i$, and
thus no vertex in $E^0_\alpha$ is the base of precisely one simple closed
path in $E_\gamma$. Continuing this process as outlined in (3.1), we find
that there is some $\beta\ge\alpha$ in $A$
such that $E_\beta$ satisfies Condition (K).

Now suppose that $E$ satisfies Condition (L). Given $\alpha\in A$, let
$p_1,p_2,\dots$ be a list of all the cycles in $E_\alpha$. (There are at
most countably many.) Each $p_i$ has an exit in $E$, say $e_i$. There is
some $\gamma\ge\alpha$ in $A$ such that $E_\gamma$ contains all the
$e_i$, and thus each cycle in $E_\alpha$ has an exit in $E_\gamma$. As
before, it follows that there is some $\beta\ge\alpha$ in $A$
such that $E_\beta$ satisfies Condition (L). \qed\enddemo

We now extend Tomforde's {\it Cuntz-Krieger Uniqueness Theorem\/} to
uncountable graphs.

\proclaim{Theorem 3.6} {\rm (Extending \cite{\Tom, Theorem 6.8,
Corollary 6.10})} Let $K$ be a field and $E$ a graph satisfying Condition
{\rm(L)}. Then every nonzero ideal of $L_K(E)$ has nonempty intersection
with $E^0$. Consequently, any ring
homomorphism $\theta: L_K(E) \rightarrow R$ satisfying $\theta(v)\ne 0$
for all
$v\in E^0$ must be injective. \endproclaim

\demo{Proof} Write $E= \dirlim E_\alpha$ as in (2.8), and set
$\Gamma= \{\gamma\in A\mid E_\gamma \text{\ satisfies Condition (L)}\}$.
In view of Lemma 3.5, $E$ is the directed union of the
subgraphs $E_\gamma$ for $\gamma\in \Gamma$, and so $L_K(E)$ is a direct
limit of the algebras $L_K(E_\gamma)$. By \cite{\Tom,
Corollary 6.10}, when $\gamma\in\Gamma$, every nonzero ideal of 
$L_K(E_\gamma)$ has nonempty intersection with $E^0_\gamma$. The theorem
follows, just as in the proof of Theorem 3.2.
\qed\enddemo

In order to adapt our modus operandi to ring-theoretic properties, we
need to express Leavitt path algebras as directed unions of countable
subrings which are themselves Leavitt path algebras. To do so, we just
combine unions of countable subgraphs with unions of countable subfields.

\definition{3.7\. Reduction to Countable Subfields} Let $\Fld$ denote the
category of fields and $\Rng$ the category of arbitrary
rings. In $\Rng$, arbitrary maps preserving addition and multiplication
are allowed as morphisms, while in $\Fld$, morphisms must be unital (in
order to qualify as field homomorphisms). Thus, $\Fld$ is a subcategory
of $\Rng$, but not a full one. Both of these categories have arbitrary
direct limits. There is an obvious functor
$$L: \Fld\times\CKGr \longrightarrow \Rng$$
such that $L(K,E)= L_K(E)$ for any field $K$ and any graph $E$. If
$(\phi,\psi): (K,E) \rightarrow (M,F)$ is a morphism in
$\Fld\times\CKGr$, then $L(\phi,\psi): L_K(E) \rightarrow L_M(F)$ is the
composition of $L_K(\phi)$ with the map $L_K(F) \rightarrow M\otimes_K
L_K(F) \equiv L_M(F)$ given by extension of scalars. We observe that the
functor $L$ preserves direct limits.

When these tools are to be used, we modify the notation of (2.8) in the
following way. Given a field $K$ and a graph $E$, write $\bigl(
(K_\alpha,E_\alpha) \bigr)_{\alpha\in A}$ for the family of ordered pairs
combining a countable subfield of $K$ with a countable CK-subgraph of
$E$. For $\alpha,\beta\in A$, define $\alpha\le\beta$ if and only if
$K_\alpha\subseteq K_\beta$ and $E_\alpha \subseteq E_\beta$. As before,
any countable ascending sequence in $A$ has a supremum.

For $\alpha\in A$, let $\lambda_\alpha$ and $\eta_\alpha$ denote the
respective inclusion maps $K_\alpha\rightarrow K$ and
$E_\alpha\rightarrow E$, and abbreviate $L(\lambda_\alpha,\eta_\alpha)$ to
$L(\eta_\alpha)$. Similarly, for $\alpha\le\beta$ in $A$, let
$\kappa_{\alpha\beta}$ and $\phi_{\alpha\beta}$ denote the respective
inclusion maps $K_\alpha\rightarrow K_\beta$ and
$E_\alpha\rightarrow E_\beta$, and abbreviate
$L(\kappa_{\alpha\beta},\phi_{\alpha\beta})$ to $L(\phi_{\alpha\beta})$.
Then  we have
$$(K,E)= \dirlim \bigl( \bigl(
(K_\alpha,E_\alpha) \bigr)_{\alpha\in A}\,,\, \bigl(
(\kappa_{\alpha\beta},\phi_{\alpha\beta}) \bigr)_{\alpha\le\beta\,
\text{in}\, A} \bigr)$$
in $\Fld\times\CKGr$, with limit maps $(\lambda_\alpha,\eta_\alpha)$, and
$$L_K(E)= \dirlim \bigl( \bigl(L_{K_\alpha}(E_\alpha) \bigr)_{\alpha\in
A}\,,\, \bigl( L(\phi_{\alpha\beta}) \bigr)_{\alpha\le\beta\,
\text{in}\, A} \bigr)$$
in $\Rng$, with limit maps $L(\eta_\alpha)$.
\enddefinition

The M.O. of (3.1) readily adapts to ring-theoretic properties in the
setting of (3.7), as in the following proof.

\proclaim{Theorem 3.8} {\rm (Extending \cite{\Tom, Theorem 6.16})}  Let
$K$ be a field and $E$ a graph. Then all ideals of $L_K(E)$ are
homogeneous if and only if $E$ satisfies Condition
{\rm(K)}. \endproclaim

\demo{Proof} $(\Longleftarrow)$: Write $E= \dirlim E_\alpha$ as in (2.8), 
and let
$\Gamma$ denote the set of those $\gamma\in A$ such that
$E_\gamma$ satisfies Condition (K). In view of Lemma 3.5, $E$ is the
directed union of the subgraphs $E_\gamma$ for $\gamma\in \Gamma$, and so
$L_K(E)$ is a direct limit of the algebras $L_K(E_\gamma)$, with limit
maps $L_K(\eta_\gamma)$. 

If $I$ is an ideal of $L_K(E)$, each $I_\gamma=
L_K(\eta_\gamma)^{-1}(I)$ is an ideal of $L_K(E_\gamma)$, and $I=
\bigcup_{\gamma\in\Gamma}  L_K(\eta_\gamma)(I_\gamma)$. By \cite{\Tom,
Theorem 6.16}, each $I_\gamma$ is a homogeneous ideal of $L_K(E_\gamma)$,
and so
$L_K(\eta_\gamma)(I_\gamma)$ is a homogeneous $K$-subspace of $L_K(E)$.
Thus their union, $I$, is a homogeneous ideal.

$(\Longrightarrow)$: This time, write $(K,E)= \dirlim
(K_\alpha,E_\alpha)$ and $L_K(E)= \dirlim L_{K_\alpha}(E_\alpha)$ as in
(3.7). Let us use the notations $\langle -\rangle_\alpha$ and
$\langle -\rangle_\infty$ for ideals in
$L_{K_\alpha}(E_\alpha)$ and $L_K(E)$, respectively.

Given $\alpha\in A$, a finite subset $X$ of $L_{K_\alpha}(E_\alpha)$, and
an element $y\in \langle X\rangle_\alpha$, we have $L(\eta_\alpha)(y)\in
\langle L(\eta_\alpha)(X)\rangle_\infty$, and so the homogeneous
components of $L(\eta_\alpha)(y)$ all lie in
$\langle L(\eta_\alpha)(X)\rangle_\infty$. At most finitely many of these
components are nonzero, and so there exists $\gamma\ge\alpha$ in $A$ such
that the homogeneous components of $L(\phi_{\alpha\gamma})(y)$ all lie in
$\langle L(\phi_{\alpha\gamma})(X) \rangle_\gamma$. Now
$L_{K_\alpha}(E_\alpha)$ has only countably many finite subsets, and each
ideal of $L_{K_\alpha}(E_\alpha)$ has only countably many elements, so in
view of the previous observations, there exists $\delta\ge\alpha$ in $A$
such that for all finite subsets $X\subseteq L_{K_\alpha}(E_\alpha)$ and
all elements $y\in \langle X\rangle_\alpha$, the homogeneous components
of $L(\phi_{\alpha\delta})(y)$ all lie in
$\langle L(\phi_{\alpha\delta})(X) \rangle_\delta$. Repeating this
procedure countably many times and taking the supremum of the resulting
indices, we obtain an index $\beta\ge\alpha$ in $A$ such that all
finitely generated ideals of $L_{K_\beta}(E_\beta)$ are homogeneous.
Consequently, all ideals of $L_{K_\beta}(E_\beta)$ are homogeneous.

Combining the above result with \cite{\Tom, Theorem 6.16}, we conclude
that for each $\alpha\in A$, there exists $\beta\ge\alpha$ in $A$ such
that $E_\beta$ satisfies Condition (K). Therefore $E$ satisfies Condition
(K). \qed\enddemo

\definition{3.9\. Saturated Hereditary Sets of Vertices} Let $E$ be a
graph and $H$ a subset (possibly empty) of $E^0$. We say that $H$ is {\it
hereditary\/} if all edges leaving $H$ end in $H$, that is, whenever
$e\in E^1$ and $s(e)\in H$, then also $r(e)\in H$. The set $H$ is called
{\it saturated\/} provided that 
\roster
\item"\bull" Each vertex $v\in E^0$ which is neither a sink nor an
infinite emitter, and which satisfies $r(e)\in H$ for all $e\in
s_E^{-1}(v)$, must lie in $H$.
\endroster
Observe that any intersection of saturated hereditary subsets of $E^0$ is
again hereditary and saturated. Hence, for any subset $V\subseteq E^0$,
we can define the {\it saturated hereditary subset of $E^0$ generated by
$V$\/} as the smallest saturated hereditary subset containing $V$. This
set can be described as follows.
\enddefinition

\proclaim{Observation 3.10} Let $E$ be a graph, $V\subseteq E^0$, and
$w\in E^0$. Then
$w$ lies in the saturated hereditary subset of $E^0$ generated by $V$ if
and only if there exists a finite sequence $w_0,w_1,\dots,w_n=w$ in
$E^0$ such that $w_0\in V$ and for all $i=1,\dots,n$, one of the following
holds:
\roster
\item $w_i\in V$;
\item There exists an edge in $E^1$ from $w_{i-1}$ to $w_i$;
\item $w_i$ is neither a sink nor an infinite emitter, and $r(e)\in
\{w_0,\dots,w_{i-1}\}$ for all $e\in s_E^{-1}(w_i)$.
\endroster
\endproclaim

Abrams and Aranda Pino characterized simplicity of Leavitt path algebras
of countable row-finite graphs in \cite{\AAa, Theorem 3.11}. They later
removed the row-finiteness hypothesis \cite{\AAc, Theorem 3.1}, as did
Tomforde, independently \cite{\Tom, Theorem 6.18}.

\proclaim{Theorem 3.11} {\rm (Extending \cite{\AAc, Theorem 3.1; \Tom,
Theorem 6.18})}  Let
$K$ be a field and $E$ a graph. Then $L_K(E)$ is a simple ring if and
only if $E$ satisfies  Condition
{\rm(L)} and the only saturated hereditary subsets of $E^0$ are
$\varnothing$ and $E^0$. \endproclaim

\demo{Proof} $(\Longleftarrow)$: Write $E= \dirlim E_\alpha$ and $L_K(E)=
\dirlim L_{K}(E_\alpha)$ as in (2.8).  Lemma 3.5 shows that for
each $\alpha\in A$, there exists
$\beta\ge\alpha$ in $A$ such that $E_\beta$ satisfies Condition (L).

If $\alpha\in A$ and $v\in E^0_\alpha$, the saturated hereditary subset
of $E^0$ generated by $v$ must equal $E^0$. In view of
Observation 3.10, we see that for each $w\in E^0_\alpha$, there exists
$\gamma\ge\alpha$ in $A$ such that $w$ lies in the
saturated hereditary subset of $E^0_\gamma$ generated by
$v$. Following the M.O. of (3.1), we conclude that
for each $\alpha\in A$, there exists $\beta\ge\alpha$ in $A$ such that
the saturated hereditary subset of $E^0_\beta$ generated by any vertex
equals $E^0_\beta$. 

Now merge the above results as follows. Given $\alpha\in A$, there is an
ascending sequence $\alpha\le \beta(1)\le \beta(2)\le \cdots$ in $A$ such
that the graphs $E_{\beta(i)}$ for $i$ even satisfy condition (L), while
those for $i$ odd satisfy the condition of the previous paragraph. The
supremum of the $\beta(i)$ is then an index $\beta\ge\alpha$ in $A$ such
that $E_\beta$ satisfies Condition (L) and the saturated hereditary
subset of $E^0_\beta$ generated by any vertex equals $E^0_\beta$. In
particular, the only saturated hereditary subsets of $E^0_\beta$ are
$\varnothing$ and $E^0_\beta$.  Combining this result with
\cite{\AAc, Theorem 3.1;
\Tom, Theorem 6.18}, we see that for each $\alpha\in A$, there exists
$\beta\ge\alpha$ in $A$ such that $L_K(E_\beta)$ is a simple ring.
Therefore $L_K(E)$ must be simple.

$(\Longrightarrow)$: Now write $(K,E)= \dirlim
(K_\alpha,E_\alpha)$ and $L_K(E)= \dirlim L_{K_\alpha}(E_\alpha)$ as in
(3.7). As in the proof of the previous theorem, we will use the notations
$\langle -\rangle_\alpha$ and
$\langle -\rangle_\infty$ for ideals in
$L_{K_\alpha}(E_\alpha)$ and $L_K(E)$, respectively.

We claim that for each $\alpha\in A$, there exists
$\beta\ge\alpha$ in $A$ such that $L_{K_\beta}(E_\beta)$ is a simple
ring. By our usual M.O., it is enough to show that for any $\alpha\in
A$ and any nonzero
$x,y\in L_{K_\alpha}(E_\alpha)$, there is some $\gamma\ge\alpha$ in $A$
such that $L(\phi_{\alpha\gamma})(x) \in \langle
L(\phi_{\alpha\gamma})(y)\rangle_\gamma$. But since $L(\eta_\alpha)$ is
injective (Corollary 3.3), the element $L(\eta_\alpha)(y)\in L_K(E)$ is
nonzero, whence
$L(\eta_\alpha)(x) \in \langle L(\eta_\alpha)(y) \rangle_\infty$. Since
$L_K(E)= \dirlim L_{K_\alpha}(E_\alpha)$, it follows
that $L(\phi_{\alpha\gamma})(x) \in \langle
L(\phi_{\alpha\gamma})(y)\rangle_\gamma$ for some $\gamma\ge\alpha$ in
$A$, as desired.

Combining the above claim with \cite{\AAc, Theorem 3.1; \Tom, Theorem
6.18}, we find that for each $\alpha\in A$, there is some
$\beta\ge\alpha$ in $A$ such that
$E_\beta$ satisfies Condition (L). It immediately follows that $E$
satisfies Condition (L).

Suppose that $E^0$ contains a proper nonempty saturated hereditary subset
$H$. Then there exist $\alpha\in A$ and $v,w\in E^0_\alpha$ such that
$\eta^0_\alpha(v)\in H$ while $\eta^0_\alpha(w)\notin H$. In view of the
claims above, we may assume that $L_{K_\alpha}(E_\alpha)$ is simple.
Since $\eta_\alpha$ is a CK-morphism, the set $H_\alpha=
(\eta^0_\alpha)^{-1}(H)$ is a saturated hereditary subset of
$E^0_\alpha$. But $v\in H_\alpha$ while $w\notin H_\alpha$, which
contradicts \cite{\AAc, Theorem 3.1; \Tom, Theorem 6.18}. Therefore $E^0$
contains no proper nonempty saturated hereditary subsets. \qed\enddemo

\head 4. Exchange Rings\endhead

\definition{4.1\. Exchange Rings} In the unital setting, exchange rings
are rings over which the regular representation (the standard free module
of rank
$1$) satisfies the exchange property in direct sum decompositions. The
definition of the exchange property and many consequences can be found in
numerous papers, of which we mention \cite{\War}, \cite{\Nic}, and the
survey
\cite{\AraB}. For present purposes, the key point is that exchange rings
can be described by finitely many ring-theoretic equations \cite{\GW, p\.
167; \Nic, Theorem 2.1}: A unital ring
$R$ is an exchange ring if and only if for each $x\in R$, there exists an
idempotent $e\in xR$ such that $1-e\in (1-x)R$. Non-unital exchange rings
were introduced by Ara \cite{\Ara}, who defined a ring $R$ to be an {\it
exchange ring\/} provided that for each $x\in R$, there exist elements
$e,r,s\in R$ such that $e$ is an idempotent and $e=xr= x+s-xs$. In
case $R$ is a ring which is generated as an ideal by its idempotents, it
follows from
\cite{\AGS, Theorem 3.3} that
$R$ is an exchange ring if and only if all the corners $eRe$, for
idempotents $e\in R$, are (unital) exchange rings.
\enddefinition

In \cite{\APS, Theorem 4.5}, Aranda Pino, Pardo, and Siles Molina showed
that the Leavitt path algebra of a countable row-finite graph $E$,
over any base field, is an exchange ring if and only if $E$ satisfies
Condition (K). The row-finiteness hypothesis was removed by Abrams and
Aranda Pino in \cite{\AAc, Theorem 5.4}. We can now remove the
countability assumption.

\proclaim{Theorem 4.2} {\rm (Extending \cite{\AAc, Theorem 5.4})} Let $K$
be a field and $E$ a graph. Then $L_K(E)$ is an exchange ring if and only
if $E$ satisfies Condition {\rm(K)}.  \endproclaim

\demo{Proof}  $(\Longleftarrow)$: Write $E= \dirlim E_\alpha$ as in (2.8), 
and let
$\Gamma$ denote the set of those $\gamma\in A$ such that
$E_\gamma$ satisfies Condition (K). In view of Lemma 3.5, $E$ is the
directed union of the subgraphs $E_\gamma$ for $\gamma\in \Gamma$, and so
$L_K(E)$ is a direct limit of the algebras $L_K(E_\gamma)$. Each
$L_K(E_\gamma)$ is an exchange ring by \cite{\AAc, Theorem 5.4}, and
therefore $L_K(E)$ must be an exchange ring.

$(\Longrightarrow)$: Now write $(K,E)= \dirlim
(K_\alpha,E_\alpha)$ and $L_K(E)= \dirlim L_{K_\alpha}(E_\alpha)$ as in
(3.7). 

Given $\alpha\in A$ and $x\in L_{K_\alpha}(E_\alpha)$, there exist
elements $e,r,s\in L_K(E)$ such that $e$ is an idempotent and
$e= L(\eta_\alpha)(x)r= L(\eta_\alpha)(x)+ s- L(\eta_\alpha)(x)s$. There
exist $\gamma\ge\alpha$ in $A$ and $e',r',s'\in L_{K_\gamma}(E_\gamma)$
such that $e'$ is an idempotent and $e'= L(\phi_{\alpha\gamma})(x)r'=
L(\phi_{\alpha\gamma})(x)+ s'- L(\phi_{\alpha\gamma})(x)s'$. By our usual
M.O., it follows that for each $\alpha\in A$, there is some
$\beta\ge\alpha$ in $A$ such that $L_{K_\beta}(E_\beta)$ is an exchange
ring. For each such $\beta$, the graph $E_\beta$ satisfies Condition (K)
by \cite{\AAc, Theorem 5.4}. Therefore $E$ satisfies condition (K). 
\qed\enddemo

\head 5. Non-Stable K-Theory \endhead

The group $K_0$ of a ring $R$, precisely because it is a group, ignores
direct sum cancellation questions. In particular, taking the unital case
for the moment, finitely generated projective $R$-modules $A$ and $B$
represent the same class in
$K_0(R)$ if and only if they are {\it stably\/} isomorphic, meaning that
$A\oplus C\cong B\oplus C$ for some finitely generated
projective $R$-module $C$. In order to keep track of isomorphism
classes, one builds a monoid, denoted $V(R)$, in place of the group
$K_0(R)$. This construction, which we now sketch, is a central object in
what has become known as ``non-stable K-theory''. 

\definition{5.1\. $V(-)$} Let $I$ be a ring (with or without unit). There
are two equivalent constructions of the monoid $V(I)$, one via projective
modules over a unital extension ring of $I$, one via idempotent matrices.
We begin with the second construction.

Write $a\oplus b$ for the block sum of square matrices $a$ and $b$, that
is, for the matrix $\left( \smallmatrix a&0\\ 0&b \endsmallmatrix
\right)$ where the $0$s are rectangular zero matrices of appropriate
sizes. We can view each matrix
ring $M_n(I)$ as a subring of $M_{n+1}(I)$ by identifying any $a\in
M_n(I)$ with
$a\oplus 0$. Set $M_\infty(I)=
\bigcup_{n=1}^\infty M_n(I)$. In the ``idempotent picture'' of $V(I)$,
the elements of $V(I)$ are equivalence classes $[e]$ of idempotents from
$M_\infty(I)$, where idempotents $e$ and $f$ are {\it equivalent\/} if
and only if there exist $a,b\in M_\infty(I)$ such that
$ab=e$ and $ba=f$. This set of equivalence classes becomes an abelian
monoid with the addition operation induced from block sums, that is,
$[e]+[f]= [e\oplus f]$. 

Functoriality of $V$ is clear from the above construction: Any morphism
$\phi: I\rightarrow J$ in $\Rng$ induces a morphism $M_\infty(\phi):
M_\infty(I) \rightarrow M_\infty(J)$ which preserves block sums,
idempotents, and equivalence, so $M_\infty(\phi)$ in turn induces a
monoid homomorphism $V(\phi): V(I) \rightarrow V(J)$. Thus, we obtain a
functor $V(-)$ from $\Rng$ to the category of abelian monoids. It is
a routine observation that this functor preserves direct limits.

For the ``projective picture'', choose any unital ring $R$ that contains
$I$ as a two-sided ideal. Let $\RMod$ denote the usual category of unital
left $R$-modules and module homomorphisms, and
$\FP(I,R)$ the full subcategory of $\RMod$ whose objects are those
finitely generated projective left $R$-modules $P$ such that $P=IP$. Then
$V(I)$ can be defined as the monoid of isomorphism classes of objects in
$\FP(I,R)$, with addition induced from direct sum. (In short, $V(I)$ is
the {\it Grothendieck monoid\/} of the category $\FP(I,R)$.) There is a
natural isomorphism from the previous incarnation of
$V(I)$ to this one, under which the equivalence class of an idempotent
$e\in M_n(I)$ is mapped to the isomorphism class of the module $R^ne$. In
particular, this shows that, up to isomorphism, the projective module
form of
$V(I)$ does not depend on the choice of unital ring $R$ in which $I$ is
embedded as an ideal.

The projective picture of $V(I)$ is convenient for dealing with Morita
equivalence, as follows.
\enddefinition

\definition{5.2\. Nonunital Morita Equivalence} As in the unital case,
Morita equivalence is based on equivalences of module categories. However,
the category of arbitrary left modules over a non-unital ring
$I$ is too large for the purpose -- for one thing, it contains all abelian
groups (viewed as $I$-modules with zero module multiplication). We follow
the common practice (see \cite{\GS}, for instance) in defining $\IMod$ to
be the category of those
left $I$-modules $M$ which are
\roster
\item"\bull" {\it full\/}: $IM=M$, and
\item"\bull" {\it nondegenerate\/}: $Ix=0$ implies $x=0$, for $x\in M$.
\endroster
(The morphisms in $\IMod$ are arbitrary module homomorphisms between the
above modules.) Observe that $\IMod$ has finite products (built as direct
products) and arbitrary coproducts (built as direct sums). 

Rings $I$ and $J$ are defined to be {\it Morita equivalent\/} provided the
categories $\IMod$ and $\JMod$ are equivalent.
\enddefinition

When $R$ and $S$ are Morita equivalent unital rings, the
monoids $V(R)$ and $V(S)$ are clearly isomorphic -- this follows easily
from the projective picture, since the Morita equivalence implies that the
categories of finitely generated projective modules over $R$ and $S$ are
equivalent. In order to make a similar argument for nonunital rings, we
need to deal with rings $I$ for which we can show that $V(I)$ is
isomorphic to the monoid of isomorphism classes of objects from some
categorically defined subcategory of $\IMod$. Idempotent rings are
suitable for this purpose, as follows.

\definition{5.3\. Idempotent Rings} A ring $I$, when viewed as a left
module over itself, might not be either full or nondegenerate, i.e., it
can fail to be an object in $\IMod$. Fullness occurs exactly when $I$ is
{\it idempotent\/}, that is, $I^2=I$. Nondegeneracy can either be
assumed, or obtained by factoring out a suitable ideal. If $I$ is
idempotent and $J= \{x\in I \mid Ix =0\}$, then $I/J$ is full and
nondegenerate as either a left $I$-module or a left $(I/J)$-module.

Idempotence by itself is already helpful in working with $\IMod$, as the
following observation shows. We thank P. Ara for communicating it to us.
\roster
\item"\bull" If $I$ is idempotent, then every epimorphism in $\IMod$ is
surjective.
\endroster
Given an epimorphism $f:M\rightarrow N$ in $\IMod$, set $X= \{x\in N\mid
Ix\subseteq f(M)\}$ and observe that $N/X$ is nondegenerate. It is also
full, because $N$ is full, whence $N/X$ is an object in $\IMod$. Now
since $f$ is an epimorphism, the quotient map $N\rightarrow N/X$ must
coincide with the zero map, and thus $N/X=0$. Consequently, $N=IN\subseteq
f(M)$, proving that $f$ is surjective.

The dual statement also holds:
\roster
\item"\bull" If $I$ is idempotent, then every monomorphism in $\IMod$ is
injective.
\endroster
Given a monomorphism $g: M\rightarrow N$ in $\IMod$, set $K=
g^{-1}(\{0\})$ (the usual kernel). Then $IK$ is a full, nondegenerate
submodule of $M$, and hence an object in $\IMod$. Now since $g$ is a 
monomorphism, the inclusion map $IK\rightarrow M$ must coincide with the
zero map, whence $IK=0$. Nondegeneracy then implies $K=0$, proving that
$g$ is injective.

As above, assume that $I$ is idempotent, and set $J= \{x\in I \mid Ix
=0\}$. Note that $J^2=0$. If $M$ is a nondegenerate left $I$-module, then
$IJM=0$ implies
$JM=0$, so that $M$ is, in a canonical way, a left $(I/J)$-module.
Consequently, the objects in $\IMod$ can be identified with the objects
in $\IJMod$, and thus we can identify these two categories. Finally, we
observe that factoring out $J$ does not harm $V(I)$:
\roster
\item"\bull" $V(I)\cong V(I/J)$.
\endroster
More precisely, $V(\pi): V(I)\rightarrow V(I/J)$ is an isomorphism, where
$\pi$ denotes the quotient map $I\rightarrow I/J$. To see that $V(\pi)$
is surjective, it suffices to lift idempotents from any matrix ring
$M_n(I/J)$ to $M_n(I)$. Since $M_n(J)^2=0$, this is a classical fact
(e.g., see \cite{\Lam, p\. 72, Proposition 1} and note that the proof
works just as well in the nonunital case). For injectivity, we need to
show that if $e$ and $f$ are any idempotents in some $M_n(I)$ whose
images are equivalent in $M_n(I/J)$, then $e$ and $f$ are equivalent.
This is also a classical fact, but we have not located a convenient
reference, so we sketch a proof. By assumption, there exist $a,b\in
M_n(I)$ such that $ab-e$ and $ba-f$ lie in $M_n(J)$. After replacing $a$
and $b$ by $eaf$ and $fbe$, we may assume that $a=eaf$ and
$b=fbe$. Now the element $ab\in eM_n(I)e$ is the sum of $e$ plus a
nilpotent, so it is a unit in that ring. Set $c=b(ab)^{-1}$, so that
$c=fce$ and $ac=e$. Since $ba$ is the sum of $f$ plus a nilpotent, it is
a unit in $fM_n(I)f$. Consequently, we can cancel the left hand $a$
factor from $aca=ea=af$ to conclude that $ca=f$. Therefore $e$ and $f$ are
equivalent, as needed.
\enddefinition

\definition{5.4\. Compact Objects} Recall
that an object $C$ in a category $\bold M$ with coproducts and a zero
object is said to be {\it compact\/} provided the following property
holds: Given any set $\bigl( X_\alpha \bigr)_{\alpha\in A}$ of objects in
$\bold M$ and any epimorphism $f: \coprod_{\alpha\in A} X_\alpha
\rightarrow C$, there exists a finite subset $B\subseteq A$ such that the
composition of
$f$ with the natural map $\coprod_{\beta\in B} X_\beta \rightarrow
\coprod_{\alpha\in A} X_\alpha$ is an epimorphism.
\enddefinition

\proclaim{Lemma 5.5} Let $I$ be an idempotent ring. Then
$V(I)$ is isomorphic to the monoid of isomorphism classes of compact
projective objects of $\IMod$ {\rm(}with addition induced from direct
sum{\rm)}. \endproclaim

\demo{Proof} Set $J= \{x\in I \mid Ix
=0\}$. As noted in (5.3), $J$ is an ideal of $I$ such that $I/J$ is
nondegenerate as a left module over itself,
$V(I)\cong V(I/J)$, and $\IMod=\IJMod$. Thus, after replacing $I$ by
$I/J$, we may assume that $I$ is a nondegenerate left $I$-module.

Let $R$ be the canonical unitification of $I$, namely the
unital ring containing $I$ as a two-sided ideal such that $R= \ZZ
\oplus I$. The forgetful functor provides a category isomorphism from
$\RMod$ to the category of arbitrary left $I$-modules \cite{\Fai,
Proposition 8.29B}, and we identify these two categories. Then
$\IMod$ is identified with the full subcategory of $\RMod$ whose objects
are the full nondegenerate left $I$-modules.

We claim that the objects of $\FP(I,R)$ are precisely the compact
projective objects of $\IMod$. Once this is proved, the lemma follows.

First, let $P$ be an object in $\FP(I,R)$. By definition, $P$ is a full
$I$-module. Further, $P$ is isomorphic to an $R$-submodule of a direct sum
of copies of
$R$, so $P=IP$ is also isomorphic to an $I$-submodule of a direct sum of
copies of $I$. Since $I$ is a nondegenerate module over
itself, $P$ is a nondegenerate $I$-module. Hence, $P$ is an object
of $\IMod$. Projectivity of $P$ in
$\RMod$, together with the surjectivity of epimorphisms in $\IMod$, now
implies that
$P$ is projective in
$\IMod$. 

Suppose that $\bigl( X_\alpha \bigr)_{\alpha\in A}$ is a set of objects
in $\IMod$ and $f: \coprod_{\alpha\in A} X_\alpha \rightarrow P$ an
epimorphism. We may view this coproduct as an internal direct sum of
modules. Since $f$ is surjective, $\sum_{\alpha\in A} f(X_\alpha)= P$.
Each $f(X_\alpha)$ is an $R$-submodule of $P$, and so ($P$ being finitely
generated) there must be a finite subset $B\subseteq A$ such that
$\sum_{\beta\in B} f(X_\beta)= P$. The composition of
$f$ with the natural map $\coprod_{\beta\in B} X_\beta \rightarrow
\coprod_{\alpha\in A} X_\alpha$ is thus an epimorphism, proving that $P$
is a compact object of $\IMod$.

Conversely, let $P$ be an arbitrary compact projective object of $\IMod$.
Since $\sum_{x\in P} Ix= IP= P$, it follows from the compactness of $P$
and the surjectivity of epimorphisms in $\IMod$ that there are finitely
many elements $x_1,\dots,x_n\in P$ such that $\sum_{j=1}^n Ix_j = P$. In
particular, $\sum_{j=1}^n Rx_j = P$, and so $P$ is finitely generated as
an $R$-module. Choose a free $R$-module $F$ and an epimorphism
$f:F\rightarrow P$ in $\RMod$. Then $IF$ is a full nondegenerate
$I$-module, and $f$ restricts to a surjective $I$-module homomorphism
$f':IF\rightarrow P$. Since $P$ is projective in $\IMod$, there is an
$I$-module homomorphism $h: P\rightarrow IF$ such that $f'h=\id_P$. But
$h$ can also be viewed as an $R$-module homomorphism $P\rightarrow F$
satisfying $fh=\id_P$, proving that $P$ is a projective $R$-module. Thus,
$P$ is an object of $\FP(I,R)$, and the claim is established. \qed\enddemo

The following corollary of Lemma 5.5 is well known among
certain researchers, but has not appeared in the literature to our
knowledge.

\proclaim{Corollary 5.6} If $I$ and $J$ are Morita equivalent idempotent
rings, then $V(I) \cong V(J)$. \endproclaim

\demo{Proof} Since $\IMod$ and $\JMod$ are equivalent, so are their
full subcategories of compact projective objects. \qed\enddemo

\definition{5.7\. Refinement Monoids} Let $V$ be an abelian monoid,
written additively. It is called a {\it refinement monoid\/} provided it
satisfies the {\it Riesz refinement property\/}: Whenever
$x_1,x_2,y_1,y_2\in V$ with $x_1+x_2= y_1+y_2$, there exist elements
$z_{ij} \in V$ for $i,j=1,2$ such that $z_{i1}+z_{i2}= x_i$ for $i=1,2$
while
$z_{1j}+z_{2j} =y_j$ for $j=1,2$. To describe some additional properties
such a monoid might enjoy, it is convenient to equip $V$ with the {\it
algebraic preorder\/} $\le$ defined by the existence of subtraction,
i.e., elements $x,y\in V$ satisfy $x\le y$ if and only if there is some
$v\in V$ such that $x+v=y$. This relation is reflexive, transitive, and
invariant under translation, the latter meaning that $x\le y$ implies
$x+z\le y+z$ for any $z\in V$.

The monoid $V$ is {\it separative\/} if it satisfies the following weak
cancellation condition: $x+x=x+y=y+y$ implies $x=y$, for any $x,y\in V$.
Equivalently, $V$ is separative if and only if $x+z=y+z$ implies $x=y$
for any $x,y,z\in V$ such that $z\le nx$ and $z\le ny$ for some $n\in
\NN$ \cite{\AGOP, Lemma 2.1}. Finally, $V$ is said to be {\it
unperforated\/} provided that $nx\le ny$ implies $x\le y$, for any
$x,y\in V$ and any $n\in\NN$.

Observe that refinement, separativity, and unperforation are all
preserved in direct limits of monoids.
\enddefinition

Ara, Moreno, and Pardo proved in \cite{\AMP} that for any row-finite
graph $E$ and any field $K$, the monoid $V(L_K(E))$ is an unperforated,
separative refinement monoid. To remove the row-finiteness restriction,
we use desingularizations to deal with countable graphs, followed by
direct limits.

\proclaim{Theorem 5.8} {\rm (Extending \cite{\AMP, Corollary 6.5})} If $K$
is a field and
$E$ a graph, then $V(L_K(E))$  is an unperforated,
separative refinement monoid. \endproclaim

\demo{Proof} Assume first that $E$ is countable. Then there exists a
desingularization $E'$ of $E$, which is a countable row-finite graph such
that $L_K(E')$ is Morita equivalent to $L_K(E)$ \cite{\AAc, Theorem 5.2}.
The desired properties hold for $L_K(E')$ by
\cite{\AMP, Corollary 6.5}, and they then transfer to $L_K(E)$ by
Corollary 5.6. 

In the general case, we have $L_K(E)= \dirlim L_K(E_\alpha)$ as in (2.8),
where the $E_\alpha$ run over the countable CK-subgraphs of $E$. By the
previous paragraph, refinement, separativity, and unperforation hold in
each $V(L_K(E_\alpha))$, and therefore they hold in $\dirlim
V(L_K(E_\alpha))$, which is isomorphic to $V(L_K(E))$. \qed\enddemo

\head Acknowledgements\endhead

The author thanks Gene Abrams, Pere Ara, and Gonzalo Aranda Pino for
useful discussions, suggestions, and correspondence.

\enddocument